\documentclass[10pt,a4paper]{article}
\usepackage[latin1]{inputenc}
\usepackage{amsmath}
\pdfoutput=1
\usepackage{amsfonts}
\usepackage{amssymb}
\usepackage{geometry}
\usepackage{mathenv}
\usepackage{color}
\usepackage{graphicx}
\usepackage[affil-it]{authblk}
\usepackage{slashbox}
\usepackage[english]{babel}
\title{Algorithm for Sensor Network Attitude Problem}
\author{M. Carmona$^1$, OJJ. Michel$^2$, J-L. Lacoume$^{1,2}$, N. Sprynski$^1$ and B. Nicolas$^2$}
\affil{$^{1}$CEA-LETI, MINATEC-Campus, 38054 Grenoble, France \\ 
$^{2}$Gipsa-lab, BP 46 F- 38402 Grenoble Cedex, France}

\begin{document}

\maketitle
\abstract{Sensor network attitude problem consists in retrieving the attitude of each sensor of a network knowing some relative orientations between pairs of sensors. The attitude of a sensor is its orientation in an absolute axis system. We present in this paper a method for solving the sensor network attitude problem using quaternion formalism which allows to apply linear algebra tools. The proposed algorithm solves the problem when all of the relative attitudes are known. A complete characterisation of the algorithm is established: spatial complexity, time complexity and robustness. Our algorithm is validated in simulations and with real experiments.}

\section{Introduction}
\label{intro}

\indent \indent  Sensor network attitude (SNA) problem consists in retrieving the attitude of each sensors of a network knowing some relative attitudes between pairs of sensors. We can make the analogy with the sensor network location (SNL) problem which is widely studied in the literature (\textcolor{blue}{\cite{Krislock}}, \textcolor{blue}{\cite{Patwari}}) and consists in retrieving sensors position from an Euclidean distance matrix eventually incomplete and noisy. \\
\indent Principal applications concerned are motion capture and vectorial waves measurement. In motion capture, the information of attitude is interesting to reconstruct the trajectory of an object or a body \textcolor{blue}{\cite{Kuipers}}. For vectorial waves, knowing the attitude allows to retrieve waves polarisation. \\
\indent A basic algorithm is an algorithm which solves the SNA problem for a complete knowledge of the relative attitudes. As for the SNL problem \textcolor{blue}{\cite{Krislock}}, having a basic algorithm allows to develop a more general algorithm, in particular when some relative attitudes are unknown. The focus of this paper is to establish a basic algorithm for the SNA problem, to characterise it and to validate it in simulation and in experiments. \\
	
In \textcolor{blue}{\textbf{section 2}}, we formalise the SNA problem using quaternion theory which allows to apply linear algebra results. \\
\indent In \textcolor{blue}{\textbf{section 3}}, we propose then an algorithm solving the SNA problem for a complete and eventually noisy relative attitudes matrix. The most important step of the method is to estimate the highest eigenvalue and an associated eigenvector of an hermitian quaternion matrix. We adapt then the classical power iteration method for complex matrices \textcolor{blue}{\cite{Watkins}} to hermitian quaternion matrices.\\
\indent In \textcolor{blue}{\textbf{section 4}}, we study time and spatial complexities of the algorithm and we prove its robustness using classical perturbation matrix results, Weyl's theorem (\textcolor{blue}{\cite{Nakatsukasa}},\textcolor{blue}{\cite{Weyl}}) and Davis-Kahan's theorem \textcolor{blue}{\cite{Davis}}. \\
\indent In \textcolor{blue}{\textbf{section 5}}, we present an experimental validation of our algorithm using inertial systems placed on a polyhedron of known geometry used as a reference.\\
\indent In \textcolor{blue}{\textbf{Appendix}}, we can find all mathematical missing details.

\section{Terminology for the SNA problem.} 

\subsection{Quaternions and Rotations.}

\indent \indent A quaternion $\mathbf{q}$ is a four components number $\mathbf{q}=q_0+q_1\mathbf{i}+q_2\mathbf{j}+q_3\mathbf{k}$ where $q_0,q_1,q_2$ and $q_3$ are real numbers and where $\mathbf{i},\mathbf{j}$ and $\mathbf{k}$ are imaginary numbers satisfying $\mathbf{i}^2=\mathbf{j}^2=\mathbf{k}^2=\mathbf{i}\mathbf{j}\mathbf{k}=-1$. The set of quaternions is denoted by $\mathbb{H}$ and can be seen as a generalisation of the complex set. The product defined on $\mathbf{i},\mathbf{j}$ and $\mathbf{k}$ induces a product on $\mathbb{H}$ generalising the complex product. With the addition, this multiplication and the multiplication with a real, $\mathbb{H}$ is a non commutative real algebra and a division ring. It is important to note that a real commutes with all of the quaternions.\\

	Let $\mathbf{q}=q_0+q_1\mathbf{i}+q_2\mathbf{j}+q_3\mathbf{k}$ be a quaternion. $q_0$ is the real part of $\mathbf{q}$. A pure quaternion is a quaternion with a null real part. We denote by $|\mathbf{q}|_2:=\left(q_0^2+q_1^2+q_2^2+q_3^2\right)^{1/2}$ the Euclidean norm of $\mathbf{q}$. $|.|_2$ is a multiplicative norm. We denote by $\overline{\mathbf{q}}:=q_0-q_1\mathbf{i}-q_2\mathbf{j}-q_3\mathbf{k}$ the conjugate of $\mathbf{q}$. We have the following properties: $\overline{\mathbf{pq}}=\overline{\mathbf{q}}  \, \overline{\mathbf{p}}$ and $\mathbf{q}\overline{\mathbf{q}}=|\mathbf{q}|_2^2$. \\
	
	A unitary quaternion is a quaternion of norm equal to 1. The set of unitary quaternions is denoted by $\mathbb{S}$. A unitary quaternion $\mathbf{q}$ can be parametrised by an angle $\theta \in [0; 2\pi[$ and a unit 3D vector $\underline{\mathbf{u}}:=[\alpha \, \beta \, \gamma]^{T} \in \mathbb{R}^3$ by $\mathbf{q}=\cos(\theta/2)+\sin(\theta/2)\mathbf{u}$ where $\mathbf{u}:=\alpha \mathbf{i}+\beta \mathbf{j}+\gamma \mathbf{k}$ is the pure and unitary quaternion deduced from $\underline{\mathbf{u}}$. This parametrisation allows to associate to every 3D rotation matrix of angle $\theta$ and vector $\underline{\mathbf{u}}$, expressed in the canonical base of $\mathbb{R}^3$, a unitary quaternion $\mathbf{q}=\cos(\theta/2)+\sin(\theta/2)\mathbf{u}$ using the following transformation \textcolor{blue}{\cite{Kuipers}}:
\begin{eqnarray}
\underline{\underline{\mathcal{R}}}(\mathbf{q})=\begin{bmatrix}
q_0^2+q_1^2-q_2^3-q_3^2 & 2\left(q_1 q_2-q_0 q_3  \right) & 2\left(q_1 q_3+q_0 q_2  \right) \\ 
2\left(q_1 q_2+q_0 q_3  \right) & q_0^2-q_1^2+q_2^3-q_3^2 & 2\left(q_2 q_3-q_0 q_1  \right) \\ 
2\left(q_1 q_3-q_0 q_2  \right) & 2\left(q_2 q_3+q_0 q_1  \right) & q_0^2-q_1^2-q_2^3+q_3^2
\end{bmatrix} 
\label{Rq}
\end{eqnarray}
For all unitary quaternions $\mathbf{q}$, $\underline{\underline{\mathcal{R}}}(\mathbf{q})$ is then a rotation matrix and $\underline{\underline{\mathcal{R}}}(\bar{\mathbf{q}})=\underline{\underline{\mathcal{R}}}(\mathbf{q})^{T}=\underline{\underline{\mathcal{R}}}(\mathbf{q})^{-1}$, where $T$ is the transposition operation. This implies that the inverse of a rotation parametrised by $\mathbf{q}$ is parametrised by $\mathbf{q}^{-1}=\bar{\mathbf{q}}$. Furthermore, for all unitary quaternions $\mathbf{p}$ and $\mathbf{q}$, $\underline{\underline{\mathcal{R}}}(\mathbf{p}\mathbf{q})=\underline{\underline{\mathcal{R}}}(\mathbf{p})\underline{\underline{\mathcal{R}}}(\mathbf{q})$. This property shows that composing 3D rotations is equivalent to multiply the associated quaternions. It is important to note that $\underline{\underline{\mathcal{R}}}$ is not an injection because $\underline{\underline{\mathcal{R}}}(\mathbf{q})=\underline{\underline{\mathcal{R}}}(-\mathbf{q})$. Indeed, we show in \textcolor{blue}{\textbf{appendix 7.1}} that there are exactly two unitary quaternions which represent the same rotation and they are opposed. Then, two quaternions could be far with respect to the Euclidean norm but they can represent the same rotation. To get over this problem, we can only deal with unitary quaternions with a positive real part. \\

	We denote by $\mathbb{H}^{M \times N}$ the set of quaternion matrix of size $M \times N$. $\mathrm{Tr}$ is the trace operator and $*$ is the transposition-conjugation operator. An hermitian matrix is a matrix $\underline{\underline{\mathbf{A}}}$ satisfying $\underline{\underline{\mathbf{A}}}=\underline{\underline{\mathbf{A}}}^*$. Because of the non commutativity of the algebra $\mathbb{H}$, we should have to consider right and left eigenvalues of every matrix. We only need to consider right eigenvalues and we will not always mention the side for further. We can show that an hermitian matrix $\underline{\underline{\mathbf{A}}} \in \mathbb{H}^{N \times N}$ has only real eigenvalues and is diagonalisable in an orthonormal base \textit{i.e.} it exists a quaternion matrix $\underline{\underline{\mathbf{U}}} \in \mathbb{H}^{N \times N}$ satisfying $\underline{\underline{\mathbf{U}}} \, \underline{\underline{\mathbf{U}}}^*=\underline{\underline{\mathbf{I}_N}}$ such that: $\underline{\underline{\mathbf{A}}}=\underline{\underline{\mathbf{U}}} \, \underline{\underline{\mathbf{\Lambda}}}  \, \underline{\underline{\mathbf{U}}}^*$ where $\underline{\underline{\mathbf{\Lambda}}}:=\mathrm{diag}\left(\lambda_1,\ldots,\lambda_N \right)$ is the diagonal matrix which contains eigenvalues of $\underline{\underline{\mathbf{A}}}$ \textcolor{blue}{\cite{Zhang}}. Finally, the Frobenius norm of $\underline{\underline{\mathbf{A}}} \in \mathbb{H}^{M \times N}$ is $||\underline{\underline{\mathbf{A}}}||_F:=\mathrm{Tr}\left(\underline{\underline{\mathbf{A}}}^* \, \underline{\underline{\mathbf{A}}} \right)^{1/2}$. 
	
\subsection{Attitude of sensors.}
	
\indent \indent	Let $\mathfrak{R}$ be the axis system of a three-components sensor and let $\mathfrak{R}_0$ be a reference axis system. The attitude of the sensor is the 3D rotation which transforms $\mathfrak{R}$ into $\mathfrak{R}_0$. Then, as a rotation can be represented by a unitary quaternion, the attitude of a sensor can also be represented by a unitary quaternion. In the sequel, we will not distinguish the attitude, the rotation and the unitary quaternion associated. \\  
	  
	We consider now $N$ three-components sensors $S_1,\ldots,S_N$ with their axis system $\mathfrak{R}_1,\ldots,\mathfrak{R}_N$ and their attitude $\mathbf{q}_1,\ldots,\mathbf{q}_N$ stored in a vector $\underline{\mathbf{Q}}=\left[\mathbf{q}_1,\ldots,\mathbf{q}_N\right]^T \in \mathbb{S}^N$ called the attitude vector. A sensor with a known attitude is called a reference. We denote by $\underline{\mathbf{Q}}_{r}$ the sub-vector of $\underline{\mathbf{Q}}$ containing quaternions associated to the references, it is the reference vector. The relative attitude between sensor $i$ and sensor $j$ is the rotation which transforms $\mathfrak{R}_i$ into $\mathfrak{R}_j$. This rotation is represented by the quaternion $\mathbf{q}_{ij}=\mathbf{q}_{i}\overline{\mathbf{q}_{j}}$. We store those quaternions in a matrix denoted by $\underline{\underline{\mathbf{O}}}$ and called the relative attitudes matrix. 

\subsection{SNA problem.}

\indent \indent	Using notions and notations defined above, the SNA problem can now be enunciated as:
\begin{center}
\textit{"How can we estimate the attitude vector $\underline{\mathbf{Q}}$ with an incomplete and noisy relative attitudes matrix $\underline{\underline{\mathbf{O}}}$ and a reference vector $\underline{\mathbf{Q}}_{r}$?"}
\end{center}
As we precise in the introduction, we only deal with a complete relative attitudes matrix in this paper. In that case, the SNA problem can be seen as an inverse problem where the direct problem is to retrieve the relation attitude matrix from the attitude vector. This is can easily be done using the following formula:
\begin{eqnarray}
\underline{\underline{\mathbf{O}}}=\underline{\mathbf{Q}} \, \underline{\mathbf{Q}}^*
\label{O=QQ}
\end{eqnarray}

\section{Basic algorithm for the SNA problem.} 	

\indent \indent	We describe in this section a basic algorithm for the SNA problem based on the quaternion theory. First, we solve the SNA problem for a non noisy case. Then, we adapted the method to take into account uncertainties on numerical computations, relative attitudes and reference attitudes. We study the performances of our algorithm in the last subsection.  

\subsection{Complete and non noisy case.}

\indent \indent	We start by proving that every vector $\underline{\mathbf{R}} \in \mathbb{S}^N$ satisfying $\underline{\mathbf{R}} \,\underline{\mathbf{R}}^{*}=\underline{\underline{\mathbf{O}}}$ is unitary-right-collinear to the attitude vector $\underline{\mathbf{Q}}$. More precisely, we have the following theorem:\\

\noindent \textbf{Theorem 1.} $\forall \underline{\mathbf{R}} \in \mathbb{S}^{N},  \underline{\mathbf{Q}} \,\underline{\mathbf{Q}}^{*}=\underline{\mathbf{R}} \,\underline{\mathbf{R}}^{*} \Leftrightarrow \exists \mathbf{s} \in \mathbb{S}, \underline{\mathbf{Q}}=\underline{\mathbf{R}}\mathbf{s}$. \\

\noindent \textit{Proof.} The proof of the "right to left" side is easy. Let $\underline{\mathbf{R}} \in \mathbb{S}^N$ satisfying assumptions of \textcolor{blue}{theorem 1} and let $\mathbf{s}=\dfrac{1}{N}\underline{\mathbf{R}}^{*} \, \underline{\mathbf{Q}} \in \mathbb{H}$. By multiplying by $\underline{\mathbf{Q}}$ to the right of the equality $\underline{\mathbf{Q}} \, \underline{\mathbf{Q}}^{*}=\underline{\mathbf{R}} \,\underline{\mathbf{R}}^{*}$ we obtain $\underline{\mathbf{Q}}=\underline{\mathbf{R}} \, \mathbf{s}$. Furthermore, this last equality implies that $\mathbf{s}$ is unitary because $\underline{\mathbf{Q}}$ and $\underline{\mathbf{R}}$ contain unitary quaternions. $\square$\\

\noindent \textcolor{blue}{Theorem 1} shows that if a particular solution of equation \textcolor{blue}{(\ref{O=QQ})} is known we can deduce the orientation vector with at least one reference. More precisely, each column of a particular solution is unitary-right-collinear to the attitude vector $\underline{\mathbf{Q}}$. Our algorithm proposes to take into account all of information contained in the relative attitudes matrix and the reference vector.\\

	We present now a method to compute a particular solution. $\underline{\underline{\mathbf{O}}}$ is an hermitian matrix and therefore it has real eigenvalues and an orthonormal right-eigenbase \textcolor{blue}{\cite{Zhang}}. Using equation \textcolor{blue}{(\ref{O=QQ})} we can note that $\underline{\underline{\mathbf{O}}}$ is of rank 1 in $\mathbb {H}$ and $\mathrm{Tr} \, \underline{\underline{\mathbf{O}}} = N$. Those remarks prove that $\underline{\underline{\mathbf{O}}}$ has two eigenvalues : $0$ of order $N-1$ and $N$ of order 1. Furthermore, for every vector $\underline{\mathbf{R}} \in \mathbb{S}^{N}$ satisfying assumptions of \textcolor{blue}{theorem 1}, we have $\underline{\underline{\mathbf{O}}} \; \underline{\mathbf{R}} = N \underline{\mathbf{R}}$. $ \underline{\mathbf{R}} $ is thus an eigenvector of matrix $\underline{\underline{\mathbf{O}}} $ associated to the only non null eigenvalue $N$. Finally, estimating a particular solution of equation \textcolor{blue}{(\ref{O=QQ})} is equivalent to estimating an eigenvector of $\underline{\underline{\mathbf{O}}} $ associated to the eigenvalue $N$ \textit{i.e.} the highest eigenvalue of $\underline{\underline{\mathbf{O}}}$. \\

	We explain now how to estimate the attitude vector $\underline{\mathbf{Q}}$ using a particular solution $\underline{\mathbf{R}}$ and the reference vector $\underline{\mathbf{Q}}_{r}$. If there is no reference, then $\underline{\mathbf{R}}$ can be taken as the solution but there is still an ambiguity due to the existence of the rotation relating $\underline{\mathbf{R}}$ and $\underline{\mathbf{Q}}$. Otherwise, let $\underline{\mathbf{R}}_{r}$ be the sub-vector of $\underline{\mathbf{R}}$ associated to the references. According to \textcolor{blue}{theorem 1}, it exists a unitary quaternion $\mathbf{s}$ such that $\underline{\mathbf{Q}}_{r}=\underline{\mathbf{R}}_{r} \, \mathbf{s}$. Then $\mathbf{s}$ can be estimated using the following equality:
\begin{eqnarray}
\mathbf{s} =\left(\underline{\mathbf{R}}_{r}^*\underline{\mathbf{R}}_{r}\right)^{-1} \left(\underline{\mathbf{R}}_{r}^*\underline{\mathbf{Q}}_{r}\right)
\label{s=LK}
\end{eqnarray}
This allows to take into account all information contained in the reference vector. The solution can finally be expressed as  $\underline{\mathbf{Q}}=\left(\underline{\mathbf{R}}_{r}^*\underline{\mathbf{R}}_{r}\right)^{-1} \underline{\mathbf{R}}\left(\underline{\mathbf{R}}_{r}^*\underline{\mathbf{Q}}_{r}\right)$.\\

	Finally, we can summarize our basic algorithm by the following list:\\
	
\noindent  \textcolor{blue}{\textbf{step 1$^*$)}}  Compute an eigenvector $\underline{\mathbf{R}}$ of $\underline{\underline{\mathbf{O}}}$ associated to $N$  \\
\noindent  \textcolor{blue}{\textbf{step 2$^*$)}}  Compute  $\mathbf{s}$\\
\noindent  \textcolor{blue}{\textbf{step 3$^*$)}}  Compute and return $\underline{\mathbf{Q}}=\underline{\mathbf{R}}\mathbf{s}$\\ 

\noindent This algorithm is theoretical. If we have to implement it, numerical uncertainties have to be taken into account. Furthermore, the relative attitudes matrix should be the issue of an estimation process, and would be thus an approximation of the theoretical relative attitudes matrix. Those remarks are also true for the reference vector. 
	
\subsection{Complete and noisy case.}

\indent \indent	 The relative attitudes matrix and the reference vector are now considered noisy and denoted by $\underline{\underline{\hat{\mathbf{O}}}}$ and $\underline{\hat{\mathbf{Q}}}_{r}$, respectively. We denote by $\underline{\hat{\mathbf{Q}}}$ the estimated attitude vector by the adapted algorithm described below. \\

	We assume that $\underline{\underline{\hat{\mathbf{O}}}}$ is still hermitian, with unitary quaternion elements and with diagonal elements equal to 1. Those assumptions are true in practice. To adapt \textcolor{blue}{\textbf{step 1$^*$}} of the theoretical algorithm, we have to note that $N$ will not be, in general, an eigenvalue of $\underline{\underline{\hat{\mathbf{O}}}}$. Then, the algorithm has to find the highest eigenvalue $\lambda_1$ of $\underline{\underline{\hat{\mathbf{O}}}}$ and an associated eigenvector $\underline{\hat{\mathbf{R}}}$. We can extract $\underline{\hat{\mathbf{R}}}_{r}$ from $\underline{\hat{\mathbf{R}}}$ and then, using $\underline{\hat{\mathbf{Q}}}_{r}$, we can compute $\hat{\mathbf{s}}$ as in \textcolor{blue}{\textbf{step 2$^*$}}:
\begin{eqnarray}
\hat{\mathbf{s}}=\left(\underline{\hat{\mathbf{R}}}_{r}^*\underline{\hat{\mathbf{R}}}_{r}\right)^{-1} \left(\underline{\hat{\mathbf{R}}}_{r}^*\underline{\hat{\mathbf{Q}}}_{r}\right)
\label{hats}
\end{eqnarray}
We can not directly return $\underline{\hat{\mathbf{Q}}}:=\underline{\hat{\mathbf{R}}} \, \hat{\mathbf{s}}$ because, in general, every component will not be a unitary quaternion. We have to normalise each component. Finally, the adapted algorithm has the following form:\\
	
\noindent \textcolor{blue}{\textbf{step 1)}}  Compute an eigenvector $\underline{\hat{\mathbf{R}}}$ of $\underline{\underline{\hat{\mathbf{O}}}}$ associated to its highest eigenvalue $\lambda_1$  \\
\noindent \textcolor{blue}{\textbf{step 2)}}  Compute $\hat{\mathbf{s}}$\\
\noindent \textcolor{blue}{\textbf{step 3)}}  Compute $\underline{\hat{\mathbf{Q}}}$\\
\noindent \textcolor{blue}{\textbf{step 4)}}  Normalise each component of $\underline{\hat{\mathbf{Q}}}=\underline{\hat{\mathbf{R}}}\hat{\mathbf{s}}$ and return it

\section{Performances.} 

\indent \indent	We study now performances of this algorithm: spatial complexity, time complexity and robustness with respect to the noise. We illustrate our result in simulations under Matlab 2008. 

\subsection{Spatial and time complexities.} 

\indent \indent First, we can note that spatial and time complexities of \textcolor{blue}{\textbf{step 2}}, \textcolor{blue}{\textbf{step 3}} and \textcolor{blue}{\textbf{step 4}} are clearly negligible in front of those of \textcolor{blue}{\textbf{step 1}}. To compute $\underline{\hat{\mathbf{R}}}$ we apply an adapted power iteration algorithm for hermitian quaternion matrices. The algorithm and the proof of its convergence are in \textcolor{blue}{appendix 7.3}. The spatial complexity is $O(N^2)$ where $N$ is the number of sensors in the network and where $O$ is the classical Landau notation for dominated functions. The time complexity is experimentally estimated to be $O(N^3)$.

\subsection{Criteria and errors associated.}

\subsubsection{Criteria.}

\indent \indent Sensor network algorithm can be seen as an optimisation problem. Indeed, it can be formulated as a minimisation of the criterion:
\begin{eqnarray}
\mathcal{C}_1(\underline{\mathbf{P}}):=||\underline{\underline{\hat{\mathbf{O}}}}-\underline{\mathbf{P}} \, \underline{\mathbf{P}}^*||_F^2, \; \underline{\mathbf{P}}  \in \mathbb{S}^{N \times 1}, \; \underline{\mathbf{P}}_r=\underline{\hat{\mathbf{Q}}}_{r}
\end{eqnarray} 
where $\underline{\mathbf{P}}_r$ is the sub-vector of $\underline{\mathbf{P}}$ associated to reference sensors. This optimisation problem has been solved for the non noisy case in \textcolor{blue}{\textbf{section 2}}. As it is difficult to solve this problem in the noisy case, we have relaxed it in \textcolor{blue}{\textbf{section 3}} by considering two steps. First, the algorithm try to minimise the criterion:
\begin{eqnarray}
\mathcal{C}_1(\underline{\mathbf{P}} )=||\underline{\underline{ \hat{\mathbf{O}}}}-\underline{\mathbf{P}} \, \underline{\mathbf{P}}^*||_F^2, \; \underline{\mathbf{P}}  \in \mathbb{H}^{N \times 1}
\label{C1}
\end{eqnarray}
The solution retained is $\underline{\mathbf{P}}:=\underline{\hat{\mathbf{R}}}$ defined as an eigenvector of $\underline{\underline{ \hat{\mathbf{O}}}}$ associated to its highest eigenvalue. Then, the algorithm minimises:
\begin{eqnarray}
\mathcal{C}_2(\mathbf{t}):=||\underline{\hat{\mathbf{Q}}}_{r} -\underline{\hat{\mathbf{R}}}_{r}\mathbf{t}||_2^2, \; \mathbf{t}  \in \mathbb{H}
\end{eqnarray}
The solution retained is $\mathbf{t}:=\hat{\mathbf{s}}$ defined by equation \textcolor{blue}{(\ref{hats})}. 

\subsubsection{Robustness and error bounds for step 1.}

\indent \indent	We derive now error bounds for variables and criteria appearing in the algorithm. For every variable $\mathbf{X}$ we denote by: 
\begin{eqnarray}
e(\mathbf{X}):=\dfrac{||\hat{\mathbf{X}}-\mathbf{X}||_F}{||\mathbf{X}||_F} 
\end{eqnarray}
the relative error obtained for the estimation of $\mathbf{X}$ by $\hat{\mathbf{X}}$. \\

	In \textcolor{blue}{\textbf{step 1}} we assume that the highest eigenvalue $\lambda_1$ corresponds to the theoretical eigenvalue $N$. This assumption is true for low level of noise. Indeed, let $\lambda_1, \ldots, \lambda_N$ be the eigenvalues of $\underline{\underline{ \hat{\mathbf{O}}}}$ indexed in the same order as eigenvalues of $\underline{\underline{ \mathbf{O}}}$ as $N, 0,\ldots,0$. To ensure that $\lambda_1$ is the "good" eigenvalue, we have to prove that the eigenvalue associated to $N$ can not be confused with others \textit{i.e.} the eigen-gap $d=\lambda_1-\lambda_2$ has to be strictly positive. Using \textcolor{blue}{Weyl's theorem}, we show in \textcolor{blue}{\textbf{appendix 7.2}} the following inequalities: 
\begin{eqnarray}
1-e(\underline{\underline{ \mathbf{O}}}) \leq \dfrac{\lambda_1}{N} \leq 1 \; ; \; \forall i \in [|2,N|], \dfrac{|\lambda_i|}{N} \leq e(\underline{\underline{ \mathbf{O}}})
\label{Weyl}
\end{eqnarray}
Those properties imply that the robustness is ensured if $e(\underline{\underline{ \mathbf{O}}}) < 1/2$. In that case $d>0$ and \textcolor{blue}{Davis-Kahan theorem}, adapted to hermitian quaternion matrix in \textcolor{blue}{\textbf{appendix 7.2}}, leads to: 
\begin{eqnarray}
e(\underline{\mathbf{R}}) \leq \frac{e(\underline{\underline{ \mathbf{O}}})}{1-2e(\underline{\underline{ \mathbf{O}}})}
\label{DK}
\end{eqnarray}
This formula ensure the robustness of the \textcolor{blue}{\textbf{step 1}}.\\

	We compute now the value of criterion $\mathcal{C}_1$ at point $\underline{\hat{\mathbf{R}}}$. Let $\underline{\mathbf{V}}_1, \ldots, \underline{\mathbf{V}}_N$ be an orthonormal eigenbase for $\underline{\underline{ \hat{\mathbf{O}}}}$ where for all $i=1 \ldots N$, $\underline{\mathbf{V}}_i$ is associated to the eigenvalue $\lambda_i$. We have:
\begin{eqnarray}
\underline{\underline{ \hat{\mathbf{O}}}}=\displaystyle{\sum_{k=1}^N} \lambda_i \underline{\mathbf{V}}_i \, \underline{\mathbf{V}}_i^*
\label{decompO}
\end{eqnarray}
We assume without loss of generality that \textcolor{blue}{\textbf{step 1}} returns $\underline{\hat{\mathbf{R}}}=\sqrt{N} \underline{\mathbf{V}}_1$. Pre-factor $\sqrt{N}$ is only chosen to obtain concise expressions for further. Then, using definition \textcolor{blue}{(\ref{C1})} and decomposition \textcolor{blue}{(\ref{decompO})}, the value of the criterion $\mathcal{C}_1$ for this vector can be expressed as:
\begin{eqnarray}
\mathcal{C}_1(\underline{\hat{\mathbf{R}}})=2N^2\left(1-\dfrac{\lambda_1}{N} \right) \leq 2N^2e(\underline{\underline{\mathbf{O}}})
\label{C1R}
\end{eqnarray}
As criterion values are absolute and not relative, the bound derived in expression \textcolor{blue}{(\ref{C1R})} depends also on the square of the number of sensors. 	

\subsubsection{Robustness for steps 2, 3 and 4.}

\indent \indent It is difficult to explicit error bounds for $e(\mathbf{s})$ and $e(\underline{\mathbf{Q}})$. However, robustness of \textcolor{blue}{\textbf{step 2}} and \textcolor{blue}{\textbf{step 3}}  are ensured because they are continuous operations. The normalisation appearing in \textcolor{blue}{\textbf{step 4}} impacts only the angle of the estimated quaternion. Indeed, we multiply the quaternion by the inverse of its norm which is real. Then, this operation does not change the direction of the associated rotation. Finally, our algorithm is robust with respect to noise on inputs $\underline{\underline{ \hat{\mathbf{O}}}}$ and $\underline{\underline{ \hat{\mathbf{K}}}}$. \\

	We can compute the value of criterion $\mathcal{C}_2$ at point $\hat{\mathbf{s}}$:
\begin{eqnarray}
\mathcal{C}_2(\hat{\mathbf{s}})=||\underline{\hat{\mathbf{Q}}}_{r}||_F^{2} -\dfrac{||\underline{\hat{\mathbf{R}}}_{r}^*\underline{\hat{\mathbf{Q}}}_{r}||_F^2}
{||\underline{\hat{\mathbf{R}}}_{r}||_F^2}
\label{C2s}
\end{eqnarray}
Then, the more $\underline{\hat{\mathbf{R}}}_{r}$ is right-collinear to $\underline{\hat{\mathbf{Q}}}_{r}$, the more the criterion 
tends to 0 and the approximation is accurate.

\subsubsection{Simulations.}

\begin{figure}[!h]
   \begin{minipage}[b]{.5\linewidth}
      \begin{center}
		\includegraphics[scale=0.5]{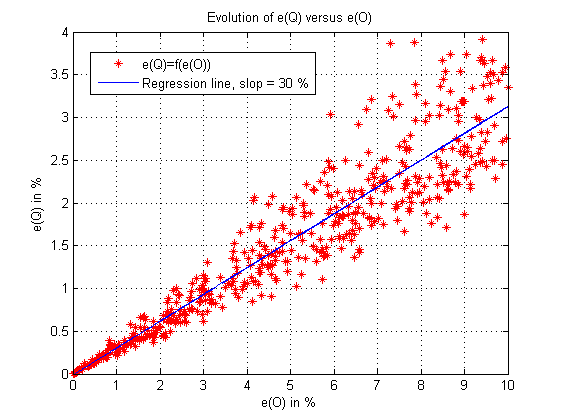}
	  \end{center}
      \caption{Estimation error.}
      \label{fig-Error}
   \end{minipage} \hfill
   \begin{minipage}[b]{.5\linewidth}
      \begin{center}
		\includegraphics[scale=0.5]{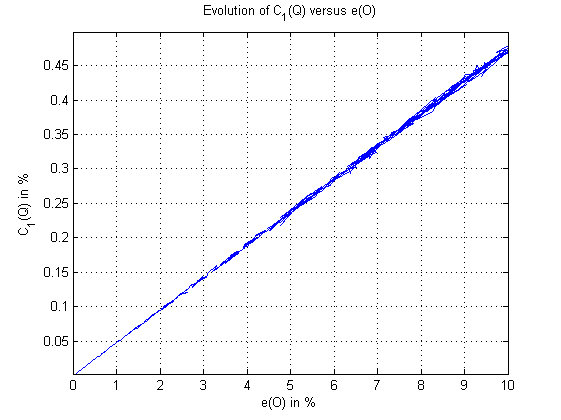}
	  \end{center}
      \caption{Criterion value versus input error.}
      \label{fig-c}
   \end{minipage}
\end{figure}

\indent \indent Inputs of the algorithm are the relative attitudes matrix $\underline{\underline{\hat{\mathbf{O}}}}$ and the reference vector $\underline{\hat{\mathbf{Q}}}_{r}$. We suppose that errors on references are negligible in front of errors on relative attitudes. This assumption is true in practice. By simulations under Matlab 2008, we observe in \textcolor{blue}{figure \ref{fig-Error}} the evolution of the output error $e(\underline{\mathbf{Q}})$ in function of the input error $e(\underline{\underline{ \mathbf{O}}})$ limited to $[0 \%,10 \%]$. Error are multiplied by 100 to be expressed in percent. We compare this evolution to a linear evolution using a linear regression. The regression line appears in \textcolor{blue}{figure \ref{fig-Error}} with a slope around of $30 \%$. Those results show that the output error is of the same order as the input error, and then validate experimentally the stability of our algorithm. We trace also  in \textcolor{blue}{figure \ref{fig-c}} the evolution of criterion value versus input error. Criterion values are divided by $N^2$ and multiply by 100 to be expressed in percent. Then, our algorithm leads to acceptable criterion values. Furthermore, that shows that the criterion is a computable indicator of estimation quality. 
	
\section{Experimental validation.}

\subsection{Experimental setup.}

\begin{figure}[!h]
\begin{center}
\includegraphics[scale=0.04]{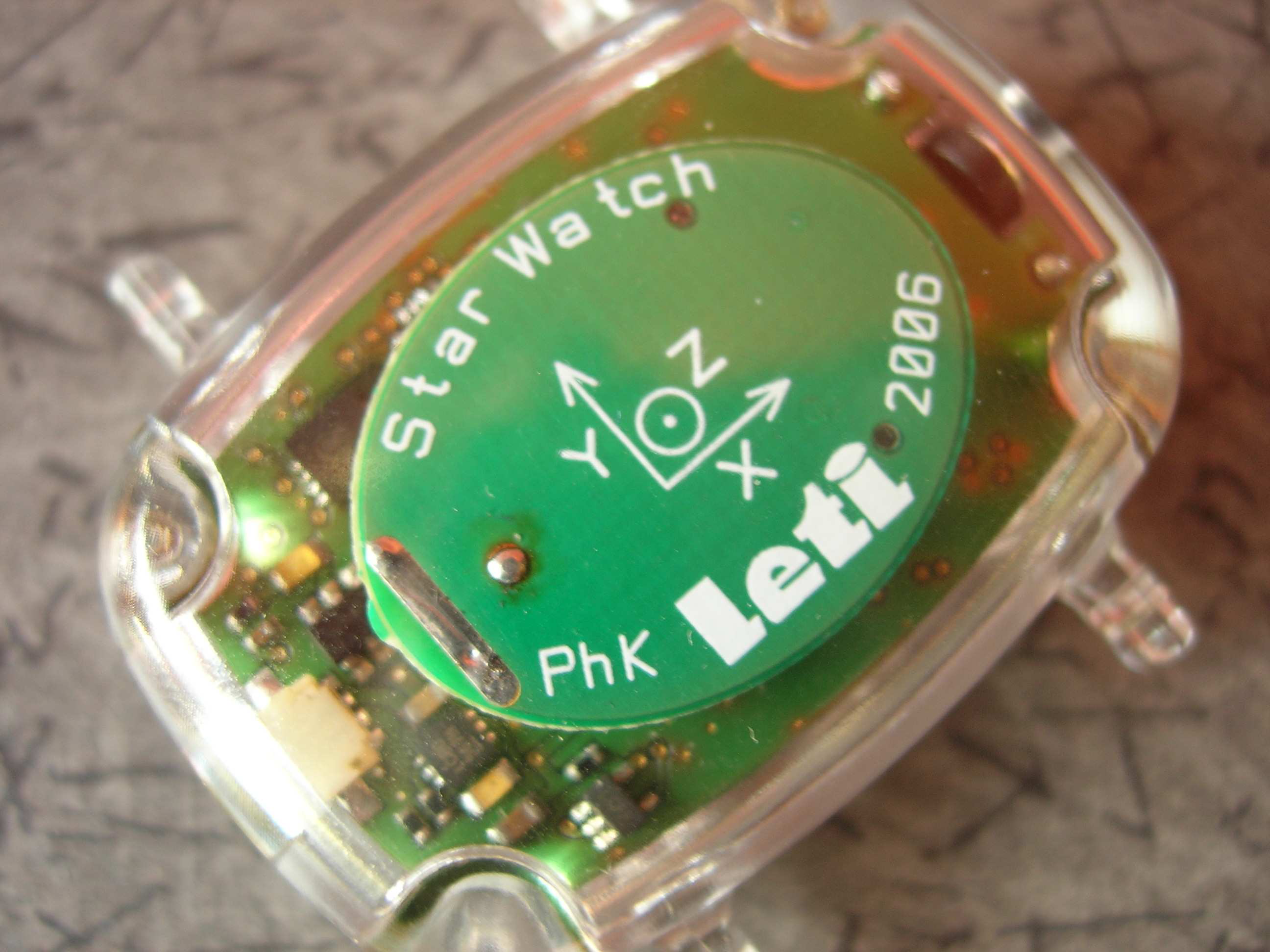} \hspace{.5cm}
\end{center}
\caption{A Star Watch.}
\label{fig-StarWatch}
\end{figure}

\indent \indent To validate experimentally our algorithm, we use systems developed by the CEA-LETI called Star Watch represented in \textcolor{blue}{figure \ref{fig-StarWatch}}. Those systems contain two sensors, a 3-components accelerometer and a 3-components magnetometer, a battery and a wireless communication module. Analogical data coming from the sensors are sampled at $200$ Hz, quantified on $12$ bits by the system it-self and sent to an acquisition system. 

\begin{figure}
\begin{center}
\includegraphics[scale=0.12]{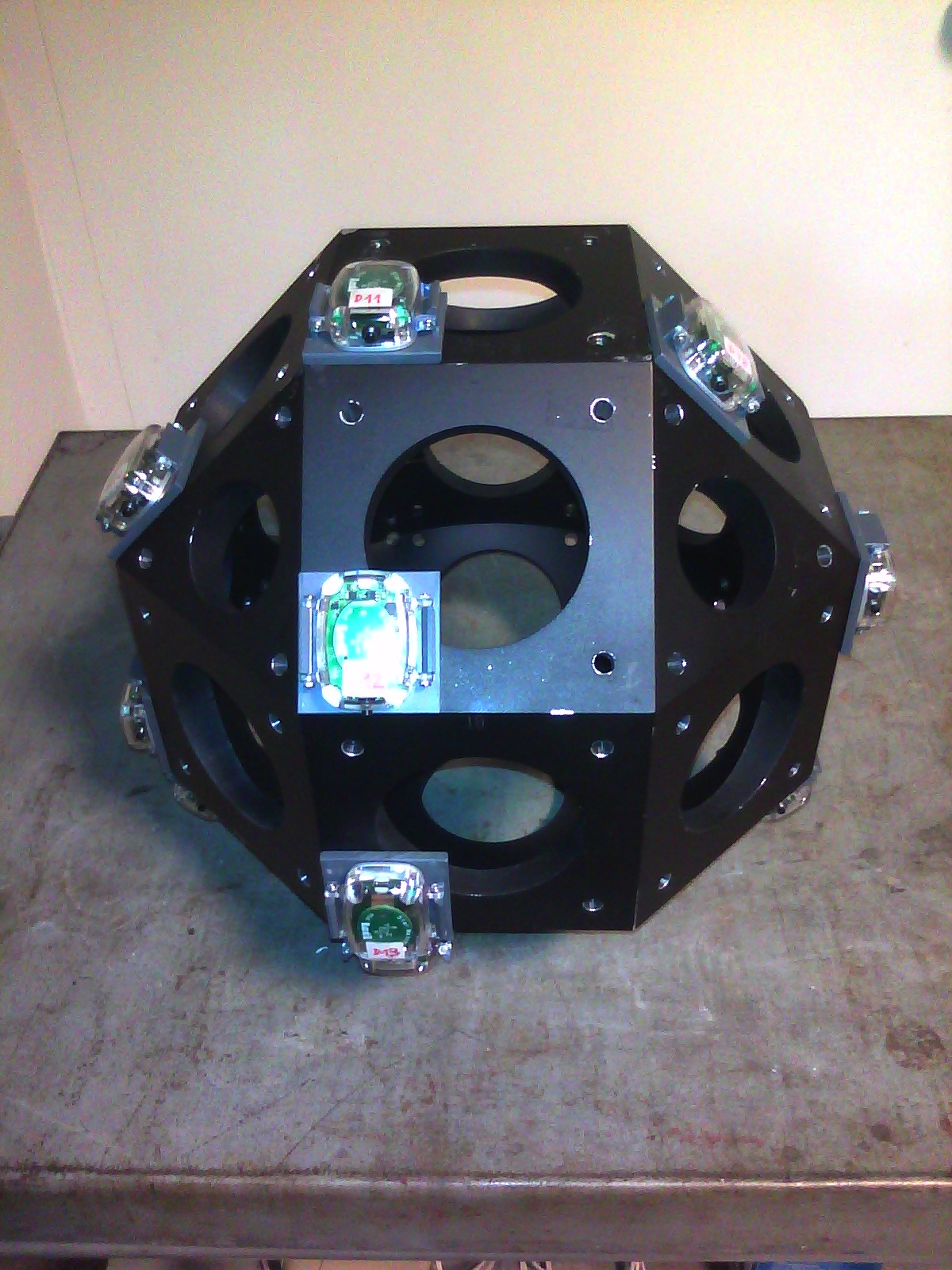}  \hspace{.5cm}
\includegraphics[scale=0.5]{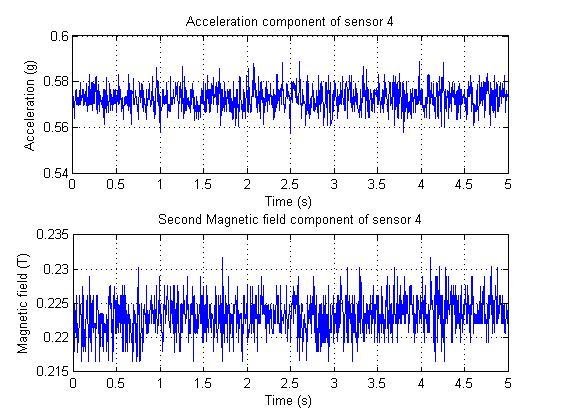}
\end{center}
\caption{\textit{Left to right}. Instrumented rhombicuboctahedron by 9 Star Watch. Measures of one component of the  accelerometer and the magnetometer of sensor 4.}
\label{fig-setup}
\end{figure}

	9 Star Watch are disposed on a rhombicuboctahedron represented in \textcolor{blue}{figure \ref{fig-setup}}. We have registered measures of the 9 static sensors during 5 seconds at 200 Hz. One component of the accelerometer and one component of the magnetometer of sensor 4 are represented in \textcolor{blue}{figure \ref{fig-setup}}. 
	
\subsection{Relative attitude estimation.}

\indent \indent	For each static sensor $i$, returned data are samples containing six measures: three components of the gravity field $[\underline{\mathbf{g}}]_{\mathfrak{R}_i}$ and three components of the magnetic terrestrial field $[\underline{\mathbf{h}}]_{\mathfrak{R}_i}$ both expressed in the axis system of the sensor $\mathfrak{R}_i$, where for all vector $\underline{\mathbf{u}} \in \mathbb{R}^3$, $[\underline{\mathbf{u}}]_{\mathfrak{R}_i}$ is the vector in $\mathbb{R}^3$ which contains coordinates of vector $\underline{\mathbf{u}}$ expressed in the axis system $\mathfrak{R}_i$. Let $\mathfrak{R}_j$ be the axis system of sensor $j$. Then, we have:
\begin{eqnarray}
[\underline{\mathbf{g}}]_{\mathfrak{R}_j}&=&\underline{\underline{\mathcal{R}}}(\mathbf{q}_{ij})[\underline{\mathbf{g}}]_{\mathfrak{R}_i} \label{eqg} \\
\left[\underline{\mathbf{h}} \right]_{\mathfrak{R}_j}&=&\underline{\underline{\mathcal{R}}}(\mathbf{q}_{ij})[\underline{\mathbf{h}}]_{\mathfrak{R}_i}  \label{eqh}
\end{eqnarray}
where $\underline{\underline{\mathcal{R}}}$ is defined in \textcolor{blue}{(\ref{Rq})} and $\mathbf{q}_{ij}$ is the unitary quaternion associated to the relative attitude between sensors $i$ and $j$. It is well-known in motion capture, that given equations \textcolor{blue}{(\ref{eqg})} and \textcolor{blue}{(\ref{eqh})}, where $[\underline{\mathbf{g}}]_{\mathfrak{R}_i}$, $[\underline{\mathbf{g}}]_{\mathfrak{R}_j}$, $[\underline{\mathbf{h}}]_{\mathfrak{R}_i}$ and $[\underline{\mathbf{h}}]_{\mathfrak{R}_j}$ are known, we can estimate $\mathbf{q}_{ij}$ \textcolor{blue}{\cite{Kuipers}}. For this estimation, we use a classical algorithm, that we call SVDQ, described in \textcolor{blue}{\cite{Kuipers}} and based on the transformation of the equations system \textcolor{blue}{(\ref{eqg})}, \textcolor{blue}{(\ref{eqh})} into a linear equations system. This new system is solved using a singular value decomposition algorithm.  

\subsection{Estimation of sensors attitude.}

\indent \indent The absolute axis system is the axis system of the upper face of the rhombicuboctahedron. Sensors attitudes are known from the geometry of the rhombicuboctahedron, we store them in the theoretical attitudes vector $\underline{\mathbf{Q}}$. From $\underline{\mathbf{Q}}$, we compute the theoretical relative attitudes matrix $\underline{\underline{\mathbf{O}}}$. It will be useful to compute the input error. \\

	For all measurements, which are stationary, we just conserve the average of the 5 seconds measurement. Then, from those measures and using SVDQ, we estimate the $36$ relative attitudes between all pairs of sensors and we stored them in the estimated relative attitudes matrix $\underline{\underline{\hat{\mathbf{O}}}}$. Using our algorithm, we compute an estimation of the attitude vector $\underline{\hat{\mathbf{Q}}}$ where sensor 1 is the unique reference.\\
	
\begin{figure}[!h]
   \begin{minipage}[b]{.5\linewidth}
      \begin{center}
		\includegraphics[scale=0.5]{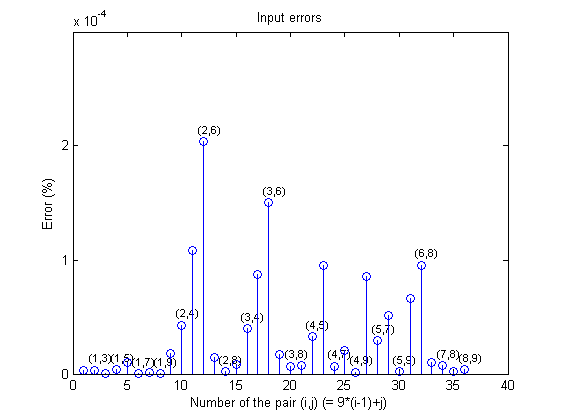}
	  \end{center}
      \caption{Input errors.}
      \label{fig-eO}
   \end{minipage} \hfill
   \begin{minipage}[b]{.5\linewidth}
      \begin{center}
		\includegraphics[scale=0.5]{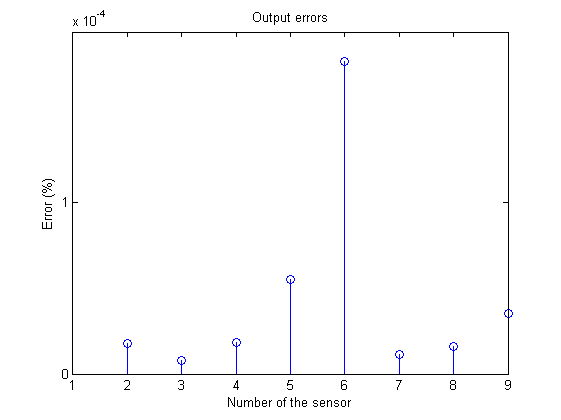}
	  \end{center}
      \caption{Output errors.}
      \label{fig-eQ}
   \end{minipage}
\end{figure}

	 According to \textcolor{blue}{\textbf{section 4}}, we compute the input error $e(\underline{\underline{\mathbf{O}}})=2 \times 10^{-4} \%$ and the output error $e(\underline{\mathbf{Q}})=1.8 \times 10^{-4} \%$. The value of the global criterion is $\mathcal{C}_1(\underline{\hat{\mathbf{Q}}})=1.7 \times 10^{-5} \%$. Those results prove the applicability of our algorithm in practical situations. For information, we trace the errors for each pair of sensors and for each sensor in \textcolor{blue}{figure \ref{fig-eO}} and \textcolor{blue}{figure \ref{fig-eQ}}. \\

We can observe in \textcolor{blue}{figure \ref{fig-eQ}} that the highest output error is on sensor 6. This is consistent with \textcolor{blue}{figure \ref{fig-eO}} where highest input errors concern pair of sensors containing this sensor. This indicates that the algorithm conserves errors repartition. 

\section{Conclusions.}

\indent \indent We presented an algorithm for solving sensors network attitude problem for a complete and eventually noisy relative attitudes matrix. The most important step relies on the estimation of the highest eigenvalue and an eigenvector associated of a quaternion hermitian matrix. We adapted then the power iteration algorithm for this type of matrix. \\
\indent Performances of the algorithm have been studied. Spatial and time complexities are $O(N^2)$ and $O(N^3)$, respectively. Robustness has been theoretically proved and validated in simulations. Simulations also give an estimation of the evolution of the output error versus the input error.  \\
\indent Finally, the algorithm has been tested in a practical situation using attitude control system (Star Watch). Those results confirm the applicability of the algorithm in real situations. \\
	
	In perspective, it will be interesting to prove the efficiency of our estimator by computing its Cramer-Rao bound when the density of probability of input noises is known. Furthermore, in order to apply this algorithm to distributed sensor networks, an important perspective is to derive a distributed version of our algorithm where the relative attitudes matrix is incomplete. 

\section{Appendix.}

\indent \indent In order to alleviate notations, we do not underline more variables when they are matrices or vectors. 

\subsection{Quaternions representing the same rotations.}
\label{sec:61}

\indent \indent To prove that it exists only two unitary quaternion which represents the same rotation, we derive the following theorem which links the error on quaternions and the error on rotation matrices associated:\\
 
\noindent \textbf{Theorem 2.} $\forall (\mathbf{q},\hat{\mathbf{q}}) \in \mathbb{S}^2$, $ ||\mathcal{R}(\hat{\mathbf{q}})-\mathcal{R}(\mathbf{q}) ||_2=|\hat{\mathbf{q}}-\mathbf{q} |_2 \sqrt{4-|\hat{\mathbf{q}}-\mathbf{q} |_2^2 }$, where $||.||_2$ is Euclidean matricial norm.  \\

\noindent \textit{Proof.} We have $||\mathcal{R}(\hat{\mathbf{q}})-\mathcal{R}(\mathbf{q}) ||_2=||\mathcal{I}_3-\mathcal{R}(\mathbf{p}) ||_2$ where $\mathcal{I}_3$ is the identity matrix of order 3 and $\mathbf{p}=\bar{\hat{\mathbf{q}}}\mathbf{q}$. $\mathbf{p}$ is a unitary quaternion, we note $\eta$ and $\underline{\mathbf{v}}$ the angle and the pure unitary quaternion associated to $\mathbf{p}$. We recall that the spectrum of a 3D rotation matrix is composed by $1, e^{\mathbf{i} \eta}, e^{-\mathbf{i} \eta}$ where $\eta$ is the angle of the rotation. As  $\mathcal{I}_3-\mathcal{R}(\mathbf{p})$ is a normal matrix, its Euclidean norm is equal to its spectral radius $\max \left \lbrace 0, |1-e^{\mathbf{i} \eta}|, |1-e^{-\mathbf{i} \eta}|\right \rbrace = 2|\sin(\eta/2)|$. We compute now $|\hat{\mathbf{q}}-\mathbf{q}|_2=|1-\mathbf{p}|_2$. As $\mathbf{p}=\cos(\eta/2)+\sin(\eta/2)\mathbf{v}$, then $|1-\mathbf{p}|_2=2|\sin(\eta/4)|$ and $\sqrt{4-|1-\mathbf{p}|_2^2}=2|\cos(\eta/4)|$. Finally, $|\hat{\mathbf{q}}-\mathbf{q} |_2 \sqrt{4-|\hat{\mathbf{q}}-\mathbf{q} |_2^2}=2|\sin(\eta/2)|=||\mathcal{R}(\hat{\mathbf{q}})-\mathcal{R}(\mathbf{q}) ||_2$. $\square$ \\

\noindent According to \textcolor{blue}{theorem 2}, its easy to show that a rotation matrix is associated only to two quaternions which are opposed. 

\subsection{Perturbation of quaternion matrices.}
\label{sec:62}

\indent \indent In order to derive properties on eigenvalues of a quaternion matrix, we use a transformation which converts a quaternion matrix into a complex matrix. Then, we can apply classical results on perturbation matrix theory. All mathematical details on quaternion matrix theory used in this appendix can be found in \textcolor{blue}{\cite{Zhang}}. \\

	A quaternion matrix $\mathbf{M}$ can be decomposed as $\mathbf{M}=\mathbf{M}_1+\mathbf{j}\mathbf{M}_2$
where $\mathbf{M}_1$ and $\mathbf{M}_2$ are complex matrices in $\mathbb{C}^{N \times N}$. This decomposition allows to associate to every matrix $\mathbf{M} \in \mathbb{H}^{N \times N}$ a complex matrix $\chi\left(\mathbf{M}\right) \in \mathbb{C}^{2N \times 2N}$ defined by:
\begin{eqnarray}
\chi(\mathbf{M}):=\begin{pmatrix}
\mathbf{M}_1&  -\overline{\mathbf{M}_2}\\ 
\mathbf{M}_2&  \overline{\mathbf{M}_1} 
\end{pmatrix} 
\end{eqnarray}
$\chi$ is $\mathbb{R}$-linear, and, satisfies $||\chi(\mathbf{M})||_F=\sqrt{2}||\mathbf{M}||_F$ and $||\chi(\mathbf{M})||_2=||\mathbf{M}||_2$ for all matrix $\mathbf{M}$. $\chi(\mathbf{M})$ is hermitian if and only if $\mathbf{M}$ is hermitian and, in that case, spectrum of $\chi(\mathbf{M})$ is equal to the spectrum of $\mathbf{M}$ where multiplicity of each eigenvalue are doubled. Furthermore, if $\mathbf{V}=\mathbf{V}_1+\mathbf{j}\mathbf{V}_2$ is an eigenvector of $\mathbf{M}$ associated to an eigenvalue $\lambda$ then $\chi(\mathbf{V}):=\left[\mathbf{V}_1 \; \mathbf{V}_2 \right]^T$  is an eigenvector of $\chi(\mathbf{M})$ associated to $\lambda$. \\

 We recall now two fundamental results of perturbation matrix theory: Weyl's theorem (\textcolor{blue}{\cite{Weyl}}, \textcolor{blue}{\cite{Nakatsukasa}}) and Davis-Kahan theorem \textcolor{blue}{\cite{Davis}} for vectorial lines:\\
 
\noindent \textbf{Theorem 3.} Let $\mathbf{A} \in \mathbb{C}^{N \times N}$ and $\hat{\mathbf{A}}=
\mathbf{A}+\delta \mathbf{A} \in \mathbb{C}^{N \times N}$ two hermitian complex matrices of eigenvalues $|a_1| >|a_2| \geq \ldots \geq |a_N|$ and $|\hat{a}_1|>|\hat{a}_2| \geq \ldots \geq |\hat{a}_N|$. Let $d=\hat{a}_1-\hat{a}_2$ be the eigengap of $\delta \mathbf{A}$. We denoted by $\mathbf{V}$ and $\hat{\mathbf{V}}$ unitary eigenvectors of $\mathbf{A}$ and $\hat{\mathbf{A}}$ associated to $a_1$ and $\hat{a}_1$, respectively. Then, we have:\\

\indent \textbf{Weyl's theorem.} $\forall k \in [|1,N|], |\hat{a}_k-a_k|\leq || \delta \mathbf{A}||_2$ \\

\indent \textbf{Davis-Kahan theorem.} If $d >0$ then $| \sin(\hat{\mathbf{V}},\mathbf{V}) | \leq \dfrac{||\delta \mathbf{A}||_F}{d}$ where $| \sin(\hat{\mathbf{V}},\mathbf{V}) |:=\dfrac{\sqrt{2}}{2}||\hat{\mathbf{V}}-\mathbf{V}||_F$.\\

\vspace{.5cm}

	As $\mathbf{O}$ and $\hat{\mathbf{O}}$ are quaternion hermitian matrices, $\chi(\mathbf{O})$ and $\chi(\hat{\mathbf{O}})$ are hermitian matrices. Eigenvalues of $\chi(\mathbf{O})$ are: $N$ of order 2 and $0$ of order $2(N-1)$. Eigenvalues of $\chi(\hat{\mathbf{O}})$ are: $\lambda_1,\ldots,\lambda_N$ where each eigenvalue is of order 2. $\chi(\mathbf{O})$ and $\chi(\hat{\mathbf{O}})$ satisfy assumptions of \textcolor{blue}{Weyl's theorem} which leads to:
\begin{eqnarray}
|\lambda_1-N| \leq ||\hat{\mathbf{O}}-\mathbf{O}||_2  \; ; \; \forall i \in [|2,N|], |\lambda_i| \leq ||\hat{\mathbf{O}}-\mathbf{O}||_2 
\label{WeylA}
\end{eqnarray}	
By diving each inequalities in \textcolor{blue}{(\ref{WeylA})} by $||\mathbf{O}||_F=N$ and by noting that $||.||_2 \leq ||.||_F$, we easily obtain inequalities in \textcolor{blue}{(\ref{Weyl})}. Inequality $|\lambda_1|\leq N$ is justified by the fact that for every matrix, its spectral radius is lower or equal to its Euclidean norm. \\

	Inequalities in \textcolor{blue}{(\ref{Weyl})}, imply that if $e(\mathbf{O}) < 1/2$ we have $d>N(1-2e(\mathbf{O}))>0$. Furthermore, we can verify that $\left|\sin(\chi(\hat{\mathbf{R}}),\chi(\mathbf{R}))\right|=e(\mathbf{R})$. With those considerations, \textcolor{blue}{Davis-Kahan theorem} applied to $\chi(\mathbf{O})$ and $\chi(\hat{\mathbf{O}})$, which satisfy the assumptions, leads to inequality  \textcolor{blue}{(\ref{DK})}. 
	
\subsection{Power iteration algorithm for a hermitian matrix.}
\label{sec:63}

\indent \indent In the complex context, power iteration method consists in estimating the highest eigenvalue $\lambda$ of a complex matrix $\mathbf{A}$ \textcolor{blue}{\cite{Watkins}}. $\lambda$ has to verify $|\lambda|>|\mu|$ for every other eigenvalue $\mu$ of $\mathbf{A}$. The algorithm gives also an eigenvector $\mathbf{V}$ associated to $\lambda$. We extend this theorem to hermitian quaternion matrix:\\

\noindent \textbf{Theorem 4.} Let $\mathbf{A} \in \mathbb{H}^{N \times N}$ be an hermitian matrix of eigenvalues $|\lambda_1| > |\lambda_2| \geq \ldots \geq |\lambda_N|$. Let $\mathbf{U}_1,\ldots,\mathbf{U}_N$ be a right-eigenbase of $\mathbb{H}^{N}$ associated to $\mathbf{A}$ . Let $\mathbf{X}_{0}=\sum_{i=1}^N \mathbf{U}_i \mathbf{a}_i $ be a vector in $\mathbb{H}^{N}$ with $\mathbf{a}_1 \neq 0$. The sequence $\mathbf{X}_{k+1}=\mathbf{A} \, \mathbf{X}_k$ satisfies $||\mathbf{X}_{k+1}||/||\mathbf{X}_k|| \to \lambda_1$ and $\mathbf{X}_{k}$ tends to an eigenvector of $\mathbf{A}$ associated to $\lambda_1$ when $k \to + \infty$. \\

\noindent \textit{Proof.} For all integers $k$, $\mathbf{X}_{k}=\mathbf{A}^k \, \mathbf{X}_0=\sum_{i=1}^N \mathbf{U}_i \lambda_i^k \mathbf{a}_i 
$, and then $\lambda_1^{-k}\mathbf{X}_{k} \to \mathbf{U}_1 \mathbf{a}_1$  when $k \to +\infty$. 
It is easy to conclude using this limit. $\square$\\

	It is interesting to notice that the power iteration method can be extended to hermitian quaternion matrices because they are diagonalisable matrices and their eigenvalues are real and then commute with every quaternion. Furthermore, even if $\mathbf{U}_1$ is unknown, the condition $\mathbf{a}_1 \neq 0$ is always satisfied due to numerical uncertainties. \\
\indent In general, to improve numerical convergence \textcolor{blue}{\cite{Watkins}}, the sequences considered are $\mathbf{V}_k=\mathbf{X}_k/||\mathbf{X}_k||$ and $\mathbf{X}_{k+1}=\mathbf{A} \, \mathbf{V}_k$ which satisfy $||\mathbf{X}_{k+1}||/||\mathbf{X}_k|| \to \lambda_1$ and $\mathbf{V}_{k}$ tends to an eigenvector of $\mathbf{A}$ associated to $\lambda_1$ when $k \to + \infty$.

\end{document}